\documentclass[11pt]{amsart}
\usepackage{amssymb, graphicx, amsthm}
\usepackage{epsfig}
\usepackage{verbatim}
\usepackage{amsfonts}
\usepackage{amsbsy}

\renewcommand{\eqref}[1]{(\ref{#1})}

\usepackage{amscd}

\begin{document}

\title[Singularities of quadratic differentials ]{Singularities of quadratic differentials and extremal Teichm\"{u}ller mappings defined by Dehn twists}
\author[C. Zhang]{C. Zhang}
\date{August 08, 2007}
\thanks{ }

\address{Department of Mathematics \\ Morehouse College
\\ Atlanta, GA 30314, USA.}
\email{czhang@morehouse.edu}

%\footnote{\footnote{}}

\maketitle

\begin{abstract}  
Let $S$ be a Riemann surface of type $(p,n)$ with $3p-3+n>0$. Let $\omega$ be a pseudo-Anosov map of $S$ that is obtained from Dehn twists along two families $\{A,B\}$ of simple closed geodesics that fill $S$. Then $\omega$ can be realized as an extremal Teichm\"{u}ller mapping on a surface of type $(p,n)$ which is also denoted by $S$. Let $\phi$ be the corresponding holomorphic quadratic differential on $S$. In this paper, we compare the locations of some distinguished points on $S$ in the $\phi$-flat metric to their locations with respect to the complete hyperbolic metric. More precisely, we show that all possible non-puncture zeros of $\phi$ must stay away from all closures of once punctured disk components of $S\backslash \{A, B\}$, and the closure of each disk component of $S\backslash \{A, B\}$ contains at most one zero of $\phi$. As a consequence of the result, we assert that the number of distinct zeros and poles of $\phi$ is less than or equal to the number of components of $S\backslash \{A, B\}$.
\end{abstract}

\bigskip
\bigskip
\bigskip

\section{Introduction}
\setcounter{equation}{0}

According to Thurston \cite{Th}, some pseudo-Anosov maps on a Riemann surface $S$ of type $(p,n)$, $3p-3+n>0$, can be constructed through Dehn twists along two simple closed geodesics (with respect to the complete hyperbolic metric on $S$ with constant curvature equal to $-1$). Let $\alpha$ and $\beta\subset S$ be two simple closed geodesics. Denote $t_{\alpha}$ and $t_{\beta}$ the positive Dehn twists along $\alpha$ and $\beta$, respectively. We assume that $\{\alpha, \beta\}$ fills $S$; that is, every component of $S\backslash \{\alpha, \beta\}$ is a disk or a once punctured disk. Thurston \cite{Th} (see also FLP \cite{FLP}) proved that the composition $t_{\alpha}\circ t_{\beta}^{-1}$ produces a pseudo-Anosov mapping class on $S$. In \cite{Lo, Pen1, Pen2} the result was generalized by showing that any word $w_0$ consisting of  $t_{\alpha}^{n_i}$ and $t_{\beta}^{-m_i}$ for positive integers $n_i, m_i$ represents a pseudo-Anosov mapping class. 

More generally, Let $A$, $B$ be families of disjoint non-trivial simple closed geodesics on $S$ so that $A\cup B$ fills $S$ in the same sense as above.  Let $w$ be any word consisting of positive Dehn twists along elements of $A$ and negative Dehn twists along elements of $B$ so that the positive Dehn twist along each element of $A$ and the negative Dehn twist along each element of $B$ occur at least once in $w$, then by a result of Penner \cite{Pen3}, $w$ also represents a pseudo-Anosov class, which means that the map $w$ can be evolved into a pseudo-Anosov map $\omega$, and if we choose $S$ properly, the map $\omega:S\rightarrow S$ is an absolutely extremal Teichm\"{u}ller mapping (see Bers \cite{Bers2}). Throughout the paper $S$ is called an $\omega$-minimal surface. Associate to $\omega$ there is a holomorphic quadratic differential $\phi$ on $S$ that may have simple poles at punctures of $S$. $\phi$ defines a flat metric on $S$. The aim of this paper is to locate (in a rather coarse manner) all possible zeros in terms of the regions obtained from cutting along the two families $\{A, B\}$ of closed geodesics  on $S$. Write 
\begin{equation}\label{COL}
S\backslash \{A,B\}=\{P_1, \ldots, P_k; Q_1, \ldots, Q_r\},
\end{equation}
\noindent where $\{Q_1,\ldots, Q_r\}$ and $\{P_1,\ldots, P_k\}$ are the collections of once punctured disk components and disk components of $S\backslash \{A,B\}$, respectively. $\{Q_1,\ldots, Q_r\}$ is empty if and only if $S$ is compact. With the notations above, the main result of this paper is the following:

\noindent {\bf Theorem 1. } {\em Let $S$ be an $\omega$-minimal surface, and $\phi$ the corresponding quadratic differential on $S$. Then: 

\noindent $(1)$ each zero of $\phi$, if exists and not a puncture of $S$, must lie in the complement of the closure of $Q_1\cup \ldots \cup Q_r$ in $S$. In particular, if $S\backslash \{A, B\}$ consists of once punctured disk components, then all zeros of $\phi$ are punctures of $S$. \\
\noindent $(2)$ the closure of each disk component $P_i$ contains at most one zero of the differential $\phi$.  }  

\noindent {\em Remark $1.1$. }By the Riemann-Roch theorem (see, for example, \cite{F-K}), if $p\geq 2$, $\phi$ has at least one zero on the compactification $\bar{S}$ of $S$. 

As a consequence of Theorem 1, we obtain:

\noindent {\bf Corollary 1. }{\em The total number of distinct zeros and poles of $\phi$ is no more than the number of the components of $S\backslash \{A, B\}$.}

\noindent {\em Remark $1.2$. } The method used to study the problem gives no information about the order of a particular zero of $\phi$. 

The idea of the proof of Theorem 1 is as follows. If a non-puncture zero $z_0$ of $\phi$ exists in $Q_1$, say, then $z_0$ would give rise to a holomorphic and isometric embedding of a Teichm\"{u}ller geodesic $\mathcal{L}\subset T(S)$ into another Teichm\"{u}ller space $T(\dot{S})$ for $\dot{S}=S\backslash \{\mbox{a point}\}$ via the fiber space $F(S)$ over $T(S)$ and a Bers isomorphism $\varphi:F(S)\rightarrow T(\dot{S})$. As a consequence, the image of $\mathcal{L}$ in $T(\dot{S})$ is a Teichm\"{u}ller geodesic in $T(\dot{S})$. A contradiction would then be derived by constructing a non-hyperbolic modular transformation that keeps invariant the image geodesic. The second statement of Theorem 1 follows from the fact that there is only one invariant geodesic under a hyperbolic transformation (I am grateful to Professor Wolpert for email communications).

\section{Preliminaries}
\setcounter{equation}{0}

We begin by reviewing some basic definitions and properties in Teichm\"{u}ller theory. Let ${\Bbb{H}}$ denote the hyperbolic plane $\{z\in {\Bbb{C}};\ \mbox{Im }z>0\}$ endowed with the hyperbolic metric
$$
ds=\frac{\left|dz\right|}{\mbox{Im }z}.
$$
\noindent  Denote $\overline{\Bbb{H}}=\{z\in {\Bbb{C}};\ \mbox{Im }z<0\}$ and let $\varrho:{\Bbb{H}}\rightarrow S$ be the universal covering with covering group $G$. Then $G$ is a torsion free finitely generated Fuchsian group of the first kind with ${\Bbb{H}}/G=S$. 
 
Let $M(G)$ be the set of Beltrami coefficients for $G$. That is, $M(G)$ consists of measurable functions $\mu$ defined on ${\Bbb{H}}$ and satisfy: (i) $\left\|\mu\right\|_{\infty}=\mbox{ess.sup }\{\left|\mu(z)\right|; \ z\in {\Bbb{H}}\}< 1$; and (ii) $\mu(g(z))\overline{g'(z)}/g'(z)=\mu(z)$ for all $g\in G$. By Ahlfors-Bers \cite{A-B}, for every $\mu\in M(G)$ there is a pair of quasiconformal maps $(w_{\mu}, w^{\mu})$ of ${\Bbb{C}}$ onto itself that satisfy these properties:

\noindent (1) both $w_{\mu}$ and $w^{\mu}$ fix $0$ and $1$; \\
(2) $\partial_{\bar{z}} w_{\mu}(z)/\partial_z w_{\mu}(z)$ is $\mu(z)$ if $z\in {\Bbb{H}}$; and $\overline{\mu(\overline{z})}$ if $z\in \overline{\Bbb{H}}$; and \\
(3) $\partial_{\bar{z}} w^{\mu}(z)/\partial_z w^{\mu}(z)$ is $\mu(z)$ if $z\in {\Bbb{H}}$ and zero if $z\in \overline{\Bbb{H}}$.

\noindent It is easy to see that $w_{\mu}$ maps ${\Bbb{H}}$ onto ${\Bbb{H}}$ while $w^{\mu}$ maps ${\Bbb{H}}$ onto a quasidisk. Two elements $\mu$ and $\nu$ in $M(G)$ are said to be equivalent if $w_{\mu}|_{\partial {\Bbb{H}}}=w_{\nu}|_{\partial {\Bbb{H}}}$, or equivalently, $w^{\mu}|_{\partial {\Bbb{H}}}=w^{\nu}|_{\partial {\Bbb{H}}}$. The equivalence class of $\mu$ is denoted by $[\mu]$. We define the Teichm\"{u}ller space $T(S)$, where $S={\Bbb{H}}/G$, is the space of equivalence classes $[\mu]$ of Beltrami coefficients $\mu\in M(G)$. $T(S)$ is a complex manifold of dimension $3p-3+n$. 

On $T(S)$ we can define the Teichm\"{u}ller distance $\left\langle [\mu],[\nu]\right\rangle$ between two points $[\mu]$ and $[\nu]$ by:
$$
\left\langle [\mu],[\nu]\right\rangle =\frac{1}{2}\ \mbox{inf}\ \left\{ \mbox{log}\ K\left(w_{\mu}\circ w_{\nu}^{-1}\right) \right\},
$$
\noindent where $K$ is the maximal dilatation of $w_{\mu}\circ w_{\nu}^{-1}$ on ${\Bbb{H}}$ and the infimum is taken through the homotopy class of $w_{\mu}\circ w_{\nu}^{-1}$ that fixes each point in $\partial {\Bbb{H}}$. The set $Q(G)$ of integrable quadratic differentials consists of holomorphic functions $\phi(z)$ on  ${\Bbb{H}}$ such that 
$$
(\phi\circ g)(z)g'(z)^2=\phi(z), \ \mbox{for all }z\in {\Bbb{H}} \  \mbox{and all } g\in G;
$$
\noindent and
$$
\left\|\phi\right\|=\frac{1}{2}\int\!\!\!\int _{\Delta}\left|\phi(z)\right|dxdy=1, 
$$ 
\noindent where $\Delta\subset {\Bbb{H}}$ is a fundamental region of $G$. Every $\phi\in Q(G)$ can be projected to a meromorphic quadratic differential on $\overline{S}$ that may have simple poles at punctures of $S$, which is denoted by $\phi$ also.  $\phi$ assigns to each uniformizing parameter $z$ a holomorphic function $\phi(z)$ such that $\phi(z)dz^2$ is invariant under a change of local coordinates. Away from zeros of $\phi$ there are naturally defined coordinates so that $\phi$ defines a flat metric that is Euclidean near every non-zero point $z$. Associate to each $\phi$ there are  horizontal and vertical trajectories defined by $\phi(z)dz^2>0$ and $\phi(z)dz^2<0$, respectively. For any $t\in (-1,1)$ and any $\phi\in Q(G)$, $t(\bar{\phi}/|\phi|)\in M(G)$. The set 
\begin{equation}\label{GEO1}
\left[t \frac{\bar{\phi}}{|\phi|}\right]\in T(S), \ \ t\in (-1,1),
\end{equation}
\noindent is a Teichm\"{u}ller geodesic. If in (\ref{GEO1}) $t$ is replaced by $z\in \Bbb{D}$, where $\Bbb{D}=\left\{ z\in {\Bbb{C}};\ |z|<1  \right\}$, we obtain a complex version of the geodesic that is also called a Teichm\"{u}ller disk. 

Notice that every self-map $\omega$ of $S$ induces a mapping class and thus a modular transformation $\chi$ that acts on $T(S)$. The collection of all such modular transformations form a group Mod$_{S}$ that is discrete and isomorphic to the group of biholomorphic maps of $T(S)$ when $S$ is not of type $(0,3),(0,4),(1,1),(1,2)$, and $(2,0)$. See Royden \cite{Ro} and Earle-Kra \cite{E-K} for more details. 

For each $\chi\in \mbox{Mod}_{S}$, Bers \cite{Bers2} introduced an index
$$
a(\chi)=\mbox{inf}_{[\mu]\in T(S)}\left\langle [\mu], \chi([\mu])\right\rangle.
$$
\noindent throughout the paper we consider those modular transformations $\chi$ for which $a(\chi)>0$. There are two cases: $a(\chi)$ is achieved and $a(\chi)$ is not achieved. In the former case, $\chi$ is called hyperbolic. In the later case, $\chi$ is called pseudo-hyperbolic. If $\chi$ is hyperbolic, then by Theorem 5 of Bers \cite{Bers2}, $a(\chi)$ assumes its value on any point in a geodesic line $\mathcal{L}$. $\chi$ keeps $\mathcal{L}$ invariant. Conversely, if an element $\chi\in \mbox{Mod}_{S}$ keeps a Teichm\"{u}ller geodesic line $\mathcal{L}$ invariant, then $\chi$ must be hyperbolic. In this case, $\chi$ is induced by a self-map of $S$, and for each Riemann surface $S$ on $\mathcal{L}$, $\chi$ is realized as an absolutely extremal self-map $\omega$ of $S$. Associate to the map $\omega$ there is an integrable meromorphic quadratic differential $\phi$ on the compactification of $S$ which is holomorphic on $S$ and may have simple poles at punctures of $S$ (see Bers \cite{Bers2}). Furthermore, $\omega$ leaves invariant both horizontal and vertical trajectories defined by $\phi$.

Topologically, the map $\omega$ is also called pseudo-Anosov that associates with a pair of transversal measured foliations determined by the quadratic differential $\phi$. By Thurston \cite{Th}, the set of pseudo-Anosov mapping classes on $S$ consists of all possible non-periodic mapping classes that do not keep any finite set of disjoint simple non-trivial closed geodesics invariant. 

The Bers fiber space $F(S)$ over $T(S)$ is the collection of pairs 
$$
\left\{\left([\mu], z \right); \ [\mu]\in T(S), \ z\in w^{\mu}({\Bbb{H}})\right\}.
$$
\noindent  The natural projection $\pi:F(S)\rightarrow T(S)$ is holomorphic. We fix a point $a\in S$ and Let $\dot{S}=S\backslash \{a\}$. Theorem 9 of \cite{Bers1} states that there is an isomorphism $\varphi:F(S)\rightarrow T(\dot{S})$.  

Let $\chi\in \mbox{Mod}_{S}$ be induced by a map $\omega:S\rightarrow S$. We lift the map $\omega$ to a map $\hat{\omega}:{\Bbb{H}}\rightarrow {\Bbb{H}}$. $\hat{\omega}$ has the property that $\hat{\omega} G {\hat{\omega}}^{-1}=G$. Suppose that $\omega':S\rightarrow S$ is another map isotopic to $\omega$. $\omega'$ can also  be lifted to a map $\hat{\omega}':{\Bbb{H}}\rightarrow {\Bbb{H}}$ that is isotopic to $\hat{\omega}$ by an isotopy fixing each point in $\partial {\Bbb{H}}$. That is, $\hat{\omega}$ and $\hat{\omega}'$ induce the same automorphism of $G$. In this case, $\hat{\omega}$ and $\hat{\omega}'$ are said to be equivalent and we denote the equivalence class of $\hat{\omega}$ by $[\hat{\omega}]$. 

Using the map $\hat{\omega}$ one constructs a biholomorphic map $\theta$ of $F(S)$ onto itself by the formula:
\begin{equation}\label{MM}
\theta ([\mu], z)=([\nu], w^{\nu}\circ \hat{\omega}\circ (w^{\mu})^{-1}(z)) \ \mbox{for every pair} \ ([\mu], z)\in F(S),
\end{equation}
\noindent  where $\nu$ is the Beltrami coefficient of $w^{\mu}\circ {\hat{\omega}}^{-1}$. 

\noindent {\bf Lemma 2.1. } (Bers \cite{Bers1}) {\em Let $\hat{\omega}$ and $\hat{\omega}'$ be in the same equivalence class. If $\hat{\omega}$ is replaced by $\hat{\omega}'$, the resulting map $\theta$ defined as $(\ref{MM})$ is unchanged. In other words, $\theta$ depends only on the equivalence class $[\hat{\omega}]$.  }\ \ \ \ \ \ \ \ \ $\Box$

%\noindent {\em Sketch of proof. } The quasidisk $w^{\mu}(\Bbb{H})$ depends only on $[\mu]$. Let $h_{\mu}$, $h_{\nu}$ be defined by the formula:
%$$
%w_{\mu}=h_{\mu}\circ w^{\mu}|_{\Bbb{H}}\ \ \mbox{and}\ \ w_{\nu}=h_{\nu}\circ w^{\nu}|_{\Bbb{H}}.
%$$
%\noindent Then $h_{\mu}:w^{\mu}({\Bbb{H}})\rightarrow {\Bbb{H}}$ and $h_{\nu}:w^{\nu}({\Bbb{H}})\rightarrow {\Bbb{H}}$ are conformal maps which depend only on $[\mu]$ and $[\nu]$, respectively. Let $\hat{\omega}$ be a lift of $\omega$. Together with (\ref{MM}), a computation shows that 
%$$
%\left([\nu], w^{\nu}\circ \hat{\omega}\circ (w^{\mu})^{-1}(z)\right)=\left([\nu], h_{\nu}^{-1}\circ A\circ h_{\mu}(z)\right),
%$$
%\noindent where $A$ is a real M\"{o}bius transformation that only depends on $[\nu]$ and the equivalence class $[\hat{\omega}]$ of $\hat{\omega}$. For the detailed argument, see \cite{Bers1}. \ \ \ $\Box$ 

Therefore, $\theta$ is uniquely determined by $[\omega]$. Lemmas 3.1 -- 3.5 of Bers \cite{Bers1} demonstrate that $\theta$ is a holomorphic automorphism of $F(S)$ that preserves the fiber structure and all such $\theta$'s form a group mod $S$ acting on $F(S)$ faithfully. 

Note that each element $g\in G$ acts on $F(S)$ by the formula:
$$
g([\mu], z)=\left([\mu], w^{\mu}\circ g\circ (w^{\mu})^{-1}(z)\right).
$$
\noindent In this way, $G$ is regarded as a normal subgroup of mod $S$, and the quotient mod $S/G$ is isomorphic to the modular group Mod$_{S}$. Let $i:\mbox{mod }S\rightarrow \mbox{Mod}_{S}$ denote the natural projection that is induced by the holomorphic projection $\pi:F(S)\rightarrow T(S)$. 

By Theorem 10 of \cite{Bers1}, the Bers isomorphism $\varphi:F(S)\rightarrow T(\dot{S})$ induces an isomorphism $\varphi^*$ of mod $S$ onto the subgroup Mod$_S^a$ of Mod$_{\dot{S}}$ that fixes the distinguished puncture $a$ via the formula:
$$
\mbox{mod } S\ni [\hat{\omega}] \ \ \stackrel{\varphi^*}{\rightarrow} \ \ \varphi\circ [\hat{\omega}]\circ \varphi^{-1}\in \mbox{Mod}_S^a. 
$$
\noindent 

\noindent The image of $[\hat{\omega}]$ in Mod$_S^a$ under $\varphi^*$ is denoted by $[\hat{\omega}]^*$.  

\section{Invariant geodesics embedded into another Teichm\"{u}ller space via a Bers isomorphism}
\setcounter{equation}{0}

Now we assume that $\chi\in \mbox{Mod}_{S}$ is a hyperbolic transformation that keeps a geodesic line $\mathcal{L}\subset T(S)$ invariant. We further assume that $S=[0]\in \mathcal{L}$. Choose $\phi\in Q(G)$ so that $\mathcal{L}=\{[t(\bar{\phi}/|\phi|)],  \  t\in (-1,1)\}$. Denote $\mu=\bar{\phi}/\phi$. Choose $\hat{\omega}:{\Bbb{H}}\rightarrow {\Bbb{H}}$ that projects to $\omega:S\rightarrow S$ that induces $\chi$. We assume that $\omega$ is an absolutely extremal Teichm\"{u}ller mapping on $S$.  By an argument of Kra \cite{Kr}, there is a hyperbolic M\"{o}bius transformation $A$ that leaves invariant $(-1,1)$ as well as ${\Bbb{D}}$ and satisfies the equation 
$$
\chi([t\mu])=\left[\mbox{Beltrami coefficient of }w^{t\mu}\circ {\hat{\omega}}^{-1}\right]=\left[A(t)\mu \right], \ \mbox{for all}\ t\in {\Bbb{D}}.
$$
\noindent Suppose that $z_0\in S$ is a zero of $\phi$ (the assumption $p\geq 2$ guarantees that $z_0\in \bar{S}$ exists by Riemann-Roch theorem, see \cite{F-K} for example).  Let $\hat{z}_0\in {\Bbb{H}}$ be such that $\varrho(\hat{z}_0)=z_0$. Let 
\begin{equation}\label{GEO2}
\hat{\mathcal{L}}=\left\{\left([t\mu], w^{t\mu}(\hat{z}_0)\right),  \ t\in (-1,1)\right\}\subset F(S).
\end{equation}  
\noindent It is easy to see that the projection $\pi:F(S)\rightarrow T(S)$ defines an embedding of $\hat{\mathcal{L}}$ into $T(S)$ with $\mathcal{L}=\pi(\hat{\mathcal{L}})$. 

The following result is well known and the argument is implicitly provided in \cite{Kr}. 

\noindent {\bf Lemma 3.1. } {\em The image $L=\varphi(\hat{\mathcal{L}})$ in $T(\dot{S})$ under the Bers isomorphism $\varphi:F(S)\rightarrow T(\dot{S})$ is a Teichm\"{u}ller geodesic.
}

\noindent {\em Proof. } By definition, $\mathcal{L}=\mathcal{L}(t)\subset T(S)$, $t\in (-1,1)$, is an isometric embedding. For any two points $x,y\in \mathcal{L}$, let $\hat{x}, \hat{y}\in \hat{\mathcal{L}}$ be such that $\pi(\hat{x})=x$ and $\pi(\hat{y})=y$. Let $x^*=\varphi(\hat{x})$ and $y^*=\varphi(\hat{y})$. By complexifying there is a Teichm\"{u}ller disk $D\subset T(S)$ with $\mathcal{L}\subset D$ and a holomorphic map $s:D\rightarrow F(S)$ defined by sending the point $[z\mu]$, $z\in {\Bbb{D}}$ and $\mu=\bar{\phi}/|\phi|$, to the point $([z\mu], w^{z\mu}(z_0))$.  It is easy to check that $\hat{\mathcal{L}}\subset s(D)$ with $s(x)=\hat{x}$ and $s(y)=\hat{y}$. Now $\varphi\circ s:D\rightarrow T(\dot{S})$ is holomorphic and is distance non-increasing. Therefore,
we have $\left\langle x^*,y^*\right\rangle\leq \left\langle x,y\right\rangle$. 

On the other hand, we notice that the natural projection $\pi:F(S)\rightarrow T(S)$ is also holomorphic, $\pi\circ \varphi^{-1}:T(\dot{S})\rightarrow T(S)$ is holomorphic and hence it is distance non-increasing. It follows that $\left\langle x,y\right\rangle\leq \left\langle x^*,y^*\right\rangle$. We thus conclude that \
$$
\left\langle x,y\right\rangle = \left\langle x^*,y^*\right\rangle,\ \mbox{for any two points}\ x,y\in \mathcal{L}.
$$ 
\noindent Hence $L=L(t)$ must also be an isometric embedding, which says that $L$ is also a Teichm\"{u}ller geodesic, as claimed. \ \ \ \ $\Box$

By \cite{Kr}, the element $\theta=[\hat{\omega}]\in \mbox{mod }S$ acts on $\hat{\mathcal{L}}$ via the formula:
$$
\theta(\hat{\mathcal{L}})=\theta\left([t\mu],\  w^{t\mu}(\hat{z}_0)\right)=\left([A(t)\mu],\  w^{A(t)\mu}\circ \hat{\omega}(\hat{z}_0)\right).
$$
\noindent From Lemma 2.1, one shows that the image $\theta(\hat{\mathcal{L}})$ in $F(S)$ only depends on $[\hat{\omega}]$. In summary, we have:

\noindent {\bf Lemma 3.2. } (Kra \cite{Kr}) {\em Suppose that $z_0\in S$ is a zero of $\phi$. An element $\theta\in \mbox{mod }S$ keeps the line $\hat{\mathcal{L}}$ invariant if and only if a representative $\hat{\omega}$ of $\theta$ can be chosen so that $\hat{\omega}(\hat{z}_0)=\hat{z}_0$.  
}

\section{Proofs of Theorem 1 and Corollary 1}
\setcounter{equation}{0}

Theorem 1 follows from Lemma 3.1 and the following result:

\noindent {\bf Theorem 2. }  {\em Let $\mathcal{L}\subset T(S)$ be a geodesic invariant under a hyperbolic transformation $\chi$. Assume that $\phi$ has a non-puncture zero $z_0$ lying in $Q_1$, say. Let $\hat{z}_0\in {\Bbb{H}}$ be such that $\varrho(\hat{z}_0)=z_0$. Then there is an element $\theta$ in mod $S$ that satisfies the following properties:

\noindent $(1)$ $\theta$ projects to $\chi$,\\
$(2)$ $\theta^*=\varphi^*(\theta)\in \mbox{Mod}_S^a$ is a non-hyperbolic but a pseudo-hyperbolic transformation acting on $T(\dot{S})$, and \\
$(3)$ the lift $\hat{\mathcal{L}}$  of $\mathcal{L}$ that passes through $\hat{z}_0$ and is defined by $(\ref{GEO2})$ is an invariant line under $\theta$. 
}

\noindent {\em Proof of Theorem $1$. } (1) If $S$ is compact, or if $\phi$ has no zeros, there is nothing to prove. So we assume that $n\geq 1$ and $z_0\in S$ is a zero of $\phi$ that is not a puncture of $S$. We further assume that $Q_1$ is a component of $S\backslash \{A,B\}$ that contains $z_0$ and a puncture $z_1$. $z_0\neq z_1$.  Let $\hat{\mathcal{L}}$ be defined in (\ref{GEO2}). By Lemma 3.1, $\varphi(\hat{\mathcal{L}})\subset T(\dot{S})$ is a Teichm\"{u}ller geodesic. By (2) and (3) of Theorem 2, $\varphi(\hat{\mathcal{L}})\subset T(\dot{S})$ is invariant under an element $\theta^*$ that is not a hyperbolic modular transformation on $T(\dot{S})$. Hence by Theorem 5 of Bers \cite{Bers2}, $\varphi(\hat{\mathcal{L}})$ is not a Teichm\"{u}ller geodesic in $T(\dot{S})$. This contradiction proves Theorem 1. 

To prove (2), we assume that there are two zeros $z_0$ and $z_1$ lying in a disk component $P_1$, say. Let $\hat{P}_1\subset {\Bbb{H}}$ be such that $\varrho|_{\hat{P}_1}:\hat{P}_1\rightarrow P_1$ is a homeomorphism. Let $\hat{z}_0, \hat{z}_1\in \hat{P}_1$ be such that $\varrho(\hat{z}_0)=z_0$ and $\varrho(\hat{z}_1)=z_1$. Let 
$\hat{\mathcal{L}}_0$ and $\hat{\mathcal{L}}_1$ be the lines passing through $\hat{z}_0$ and $\hat{z}_1$, respectively. 
From Lemma 3.1, $\varphi(\hat{\mathcal{L}}_0)$ and $\varphi(\hat{\mathcal{L}}_1)$ are distinct Teichm\"{u}ller geodesics. 

On the other hand,  by a similar argument of Theorem 2, there is a common modular transformation  $\Theta^*$ on $T(\dot{S})$ that leaves both lines $\varphi(\hat{\mathcal{L}}_0)$ and $\varphi(\hat{\mathcal{L}}_1)\subset T(\dot{S})$ invariant. By Theorem 5 of Bers \cite{Bers2}, $\Theta^*$ is hyperbolic. 
This contradicts to the uniqueness of the invariant geodesic of a hyperbolic transformation. This proves (2). \ \ $\Box$

\noindent {\em Proof of Corollary $1$. } Note that the closure $\bar{P}_i$ or $\bar{Q}_j$ of each component $P_i$ or $Q_j$ in (\ref{COL}) is a polygon with geodesic boundary segments (with respect to the hyperbolic metric). By the argument of Theorem 2, we see that each $\bar{P}_i$ or $\bar{Q}_j$ cannot contain more than one zero of $\phi$. Each pole of $\phi$ must be the puncture of some $Q_j$. By (1) of Theorem 1, each $\bar{Q}_j$ contains at most one zero $z_j$ that is the puncture of $Q_j$. In this case, $z_j$ cannot be a pole of $\phi$. Moreover, if there are two components, $\bar{P}_i$ and $\bar{P}_i$, say, in (\ref{COL}) so that a zero of $\phi$ lies in $\bar{P}_i\cap \bar{P}_j$, then the argument also shows that there do not exist any other zeros in either $\bar{P}_i$ or $\bar{P}_j$. If a zero lies in one of the intersections of $\alpha_i$ and $\beta_j$, then the same argument also shows that the closure of the related four polygons in (\ref{COL}) do not include any other zeros. Overall, we conclude that the total number of zeros and poles of $\phi$ is no more than the number of components of $S\backslash \{A,B\}$. \ \ \ \ \ \ \ \ \ $\Box$ 
 
The rest of the paper is devoted to the proof of Theorem 2. 

\section{Reducible maps projecting to pseudo-Anosov maps as a puncture is filled in}
\setcounter{equation}{0}

Let $\dot{S}$ be a Riemann surface defined in Section 2. $\dot{S}$ has type $(p,n+1)$. Let $W$ be a non periodic  non pseudo-Anosov self map of $\dot{S}$. By Thurston \cite{Th}, there exists an admissible loop system
\begin{equation}\label{LOOP}
\{c_1, c_2,\ldots, c_s\},\ s\geq 1
\end{equation}
\noindent  of simple non-trivial loops on $\dot{S}$ so that for every $i$, $1\leq i\leq s$, $W(c_i)$ is homotopic to $c_j$ for some $j$. Here by ``admissible'' we mean that no loop in (\ref{LOOP}) bounds a 1-punctured disk and  $c_i$ is not homotopic to $c_j$ whenever $i\neq j$, $1\leq i,j\leq m$. Note that $W$ may permute the components $\{R_1,\ldots R_k\}$ of $S\backslash \{c_1, c_2,\ldots, c_s\}$, and if $W$ keeps a component $R_j$ invariant, the restriction $W|_{R_j}$ could be either the identity, or periodic, or pseudo-Anosov. Thus there is an integer $N$ so that $W^N$ keeps every $c_i$ and every $R_j$ invariant, and for each $j$, $1\leq j\leq k$, $W^N|_{R_j}$ is either the identity or pseudo-Anosov. The mapping classes induced by such kind of maps are called pure mapping classes in the literature. If all $W^N|_{R_j}$, $1\leq j\leq k$, are the identity,  $W^N$ is a product of powers of positive and negative Dehn twists along certain loops in (\ref{LOOP}). In general case, $W^N$ induces a pseudo-hyperbolic transformation on $T(\dot{S})$. See Bers \cite{Bers2} for details. 

Now we consider a special case. Let $\theta=[\hat{\omega}]$ be an element of $\mbox{mod }{S}$ that projects to $\chi\in \mbox{Mod}_{S}$. We assume that $\chi$ is induced by $\omega:S\rightarrow S$ that is an absolutely extremal Teichm\"{u}ller mapping. Let $\phi\in Q(G)$ be the corresponding quadratic differential. By Royden's theorem \cite{Ro} (see also Earle-Kra \cite{E-K}), $\theta^*$ is a modular transformation. Thus $\theta^*$ is induced by a quasiconformal self-map $W$ of $\dot{S}$. The map $W$ is isotopic to $\omega$ if $W$ is viewed as a self-map of $S$. Notice that $W$ is non periodic; it may or may not be pseudo-Anosov. On the other hand, even if $W$ is pseudo-Anosov, $\dot{S}$ may not be the right candidate that makes $W$ an absolutely extremal self-map on $\dot{S}$. 

\noindent {\bf Lemma 5.1. } {\em Assume that $S$ is not compact, and $W$ is not pseudo-Anosov, then $W$ is reduced by a single closed geodesic $c_1$ that is a boundary of a twice punctured disk $\Omega\subset \dot{S}$ that encloses $a$. More precisely, if we write $\dot{S}\backslash c_1=D\cup R$, then $W|_{\Omega}$ is the identity and $W|_{R}$ is pseudo-Anosov that is essentially the same as $\omega$, and $W$ induces a pseudo-hyperbolic transformation on $T(\dot{S})$ in the sense of Bers \cite{Bers2}. }

\noindent {\em Proof. } Let $W$ be reduced by (\ref{LOOP}), and let $\gamma_i$ denote the geodesic on $S$ obtained from $c_i$ by adding the puncture $a$. Since $W$ is isotopic to $\omega$ on $S$, $\omega$ keeps the curve system $\{\gamma_1,\ldots \gamma_{s_0}\}$ invariant, where $s_0=s$ if neither two elements $c_i$ and $c_j$ bound an $a$-punctured cylinder, nor does an element $c_i$ project to a trivial loop; $s_0=s-1$ otherwise. 

Since $\omega$ is pseudo-Anosov, the set $\{\gamma_1,\ldots \gamma_{s_0}\}$ is empty. Hence the only possibility is that all geodesics in (\ref{LOOP}) are boundaries of twice punctured disks enclosing $a$. Since geodesics in (\ref{LOOP}), if not empty,  are disjoint. We must have that $s=1$ and $c_1$ in (\ref{LOOP}) is the boundary of a twice punctured disk. 

As $a$ is filled in, $c_1$ becomes a trivial loop. This means that $W|_{R}$ is essentially the same as $\omega$. Notice that $\Omega$ is a twice punctured disk, and $W|_{\Omega}$ fixes each boundary component. It follows that $W|_{\Omega}$ is isotopic to the identity. The lemma is proved. \ \ \ \ \ \ \ $\Box$

The following result along with Lemma 5.1 establishes the relationship between elements in mod $S$ and non pseudo-Anosov elements in Mod$_S^a$ via the Bers isomorphism $\varphi^*$. Recall that $[\hat{\omega}]^*=\theta^*\in \mbox{Mod}_S^a$ is induced by $W:\dot{S}\rightarrow \dot{S}$.  

\noindent {\bf Lemma 5.2. } {\em  Suppose that $S$ is not compact. If $\omega:S\rightarrow S$ does not fix any punctures of $S$, then all $W$'s must also be pseudo-Anosov. Otherwise, we assume that
$\omega:S\rightarrow S$ is pseudo-Anosov that fixes at least one puncture of $S$. Then certain non pseudo-Anosov maps $W$ of $\dot{S}$ exist, and all possible non pseudo-Anosov maps $W$ are obtained from those $\hat{\omega}:{\Bbb{H}}\rightarrow {\Bbb{H}}$ that fix a fixed point of a parabolic element $T$ of $G$. }

\noindent {\em Proof. } Assume that $\hat{\omega}:{\Bbb{H}}\rightarrow {\Bbb{H}}$ fixes the fixed point $x$ of a parabolic element $T\in G$. This implies that $\hat{\omega}\circ T=T^k\circ \hat{\omega}$ for some $k$. That is, 
\begin{equation}\label{RR}
[\hat{\omega}]^*\circ T^*={T^*}^k\circ [\hat{\omega}]^*.
\end{equation}
\noindent  From Theorem 2 of \cite{Kr, Na}, $T^*$ is represented by a Dehn twist $t_{\partial \Omega}$ along the boundary $\partial \Omega$ of a twice punctured disk $\Omega\subset \dot{S}$. Let $W:S\rightarrow S$ be a map that induces $[\hat{\omega}]^*$. From (\ref{RR}) we have 
$$
t_{W(\partial \Omega)}=t_{\partial \Omega}^k.
$$
\noindent It follows that $k=1$ and the map $W$ leaves $\partial \Omega$ invariant. So $W$ is not pseudo-Anosov. 

Conversely, assume that $W$ is not pseudo-Anosov. By Lemma 5.1, $W$ is reduced by a single geodesic $c$ that is a boundary of a twice punctured disk. This means that $W\circ t_c=t_c\circ W$. By Theorem 2 of \cite{Kr, Na} again, there is a parabolic element $T\in G$ so that $T^*=t_c$. Hence we have $[\hat{\omega}]^*\circ T^*=T^*\circ [\hat{\omega}]^*$. Thus $\hat{\omega}\circ T^k=T^k\circ \hat{\omega}$ for any integer $k$. It follows that $\hat{\omega}$ fixes the fixed point of $T$. 

In particular, if $\omega:S\rightarrow S$ does not fix any punctures of $S$, then $W$ must be pseudo-Anosov. \ \ \ \ \ \ \ $\Box$    

\noindent {\em Remark $5.1$. } The disk $\Omega\subset \dot{S}$ obtained from  Lemma $5.1$ contains another puncture $b\neq a$, which is viewed as a puncture of $S$ corresponding to the conjugacy class of $T$. Conversely, every $[\hat{\omega}]\in \mbox{mod }S$ that fixes a parabolic fixed point of $G$ produces a non pseudo-Anosov map $W$ on $\dot{S}$ that is characterized in Lemma 5.1. 

\section{pseudo-Anosov maps defined by simple closed geodesics and their lifts}
\setcounter{equation}{0}

Denote $A=\{\alpha_1,\ldots, \alpha_q\}$ and $B=\{\beta_1,\ldots,\beta_r\}$. Let $z\in S\backslash \{A,B\}$. Fix a fundamental region $\Delta$ of $G$ and let $\hat{z}=\{\varrho^{-1}(z)\}\cap \Delta$. 

Let $\hat{\alpha}_1\subset {\Bbb{H}}$ be a geodesic so that $\varrho(\hat{\alpha}_1)=\alpha_1$ and $\Delta\cap \hat{\alpha}_1\neq \oslash$ for a fixed fundamental region $\Delta$ of $G$. Note that there may be more than one choices for such a geodesic $\hat{\alpha}_1$. The geodesic $\hat{\alpha}_1$ is invariant under a simple hyperbolic element $g_{\hat{\alpha}_1}$ of $G$. Let $D_{\hat{\alpha}_1}$ and $D_{\hat{\alpha}_1}'$ be the components of ${\Bbb{H}}\backslash \hat{\alpha}_1$. We assume that $\hat{z}\in D_{\hat{\alpha}_1}'$. 

To obtain a lift $\tau_{\hat{\alpha}_1}$ of $t_{\alpha_1}$, we take an earthquake shifting along $\hat{\alpha}_1$ in such a way that it is the identity on $D_{\hat{\alpha}_1}'$; and is $g_{\hat{\alpha}_1}$ on $D_{\hat{\alpha}_1}$ away from a small neighborhood of $\hat{\alpha}_1$. We thus define $\tau_{\hat{\alpha}_1}$ everywhere on ${\Bbb{H}}$ via $G$-invariance. 

The construction of $\tau_{\hat{\alpha}_1}$ gives rise to a collection $E_{\hat{\alpha}_1}$ of half-planes, among which a partial order can be naturally defined. There are infinitely many disjoint maximal elements of $E_{\hat{\alpha}_1}$ and for each maximal element $D_{\hat{\alpha}_1}=D_{\hat{\alpha}_1}^1$ of $E_{\hat{\alpha}_1}$, there are infinitely many second level elements $D_{\hat{\alpha}_1}^2\subset D_{\hat{\alpha}_1}^1$; and for each such $D_{\hat{\alpha}_1}^2$, there are infinitely many third level elements $D_{\hat{\alpha}_1}^3$ in $D_{\hat{\alpha}_1}^2$, and so on. 

The quasiconformal homeomorphism $\tau_{\hat{\alpha}_1}$ restricts to the identity on the complement of disjoint union of all maximal elements of $\tau_{\hat{\alpha}_1}$ in ${\Bbb{H}}$; it is quasiconformal whose Beltrami coefficient is supported on (disjoint) neighborhoods of $\hat{\alpha}_1$ and its $G$-translations. Moreover, from the construction $\tau_{\hat{\alpha}_1}(z)=z$ for points $z$ on the boundaries of all maximal elements.  $\tau_{\hat{\alpha}_1}$ naturally extends to a quasisymmetric mapping of $\partial {\Bbb{H}}$ onto $\partial {\Bbb{H}}$ that fixes infinitely many hyperbolic fixed points of $G$ and infinitely many parabolic fixed points as well if $S$ is not compact. 

Let $\hat{\alpha}_1, \ldots, \hat{\alpha}_q\subset {\Bbb{H}}$ be the geodesics so that they meet in $\Delta$ in the sense that $\Delta\cap \hat{\alpha}_i\neq \oslash$ for $i=1,\ldots, q$. Since 
$\alpha_1,\ldots, \alpha_q$ are pairwise disjoint, $\hat{\alpha}_1, \ldots, \hat{\alpha}_q$ are pairwise disjoint as well.  The maximal elements $D_{\hat{\alpha}_1}, \ldots, D_{\hat{\alpha}_q}$ can be properly chosen so that 
\begin{equation}\label{ALPHA}
\hat{z}\in \Delta\backslash \{D_{\hat{\alpha}_1}, D_{\hat{\alpha}_2},\ldots, D_{\hat{\alpha}_q}\}.
\end{equation}
\noindent Since $\alpha_1, \ldots, \alpha_q$ are pairwise disjoint, $\tau_{\hat{\alpha}_i}$  commutes with $\tau_{\hat{\alpha}_j}$ for $i,j=1,\ldots,q$. Now for a positive integer tuple $\sigma=(n_1,\ldots, n_q)$, we define
$$
\hat{T}_{A}^{\sigma}=\tau_{\hat{\alpha}_1}^{n_1}\circ \tau_{\hat{\alpha}_2}^{n_2}\circ \cdots \circ \tau_{\hat{\alpha}_q}^{n_q}.
$$
\noindent We see that $\hat{T}_{A}$ does not depend on the order of those $\tau_{\hat{\alpha}_i}$. 

Similarly, let $\hat{\beta}_1, \ldots, \hat{\beta}_r\subset {\Bbb{H}}$ be the geodesics so that they meet in $\Delta$. $D_{\hat{\beta}_1}, \ldots, D_{\hat{\beta}_r}$ can also be properly chosen so that 
\begin{equation}\label{BETA}
\hat{z}\in \Delta\backslash \{D_{\hat{\beta}_1}, D_{\hat{\beta}_2},\ldots, D_{\hat{\beta}_r}\}.
\end{equation}
\noindent For a positive integer tuple $\lambda=(m_1,\ldots, m_r)$, we define
$$
\hat{T}_{B}^{\lambda}=\tau_{\hat{\beta}_1}^{m_1}\circ \tau_{\hat{\beta}_2}^{m_2}\circ \cdots \circ \tau_{\hat{\beta}_r}^{m_r}.
$$

Now we consider the map 
\begin{equation}\label{MAP}
\hat{T}_{AB}=\prod_i\left( \hat{T}_{A}^{\sigma_i}\circ \hat{T}_{B}^{-\lambda_i} \right)
\end{equation}
\noindent for positive integer tuples $\sigma_i$ and $\lambda_i$. We claim that $\hat{T}_{AB}$ defined as (\ref{MAP}) is a lift of the pseudo-Anosov mapping class: 
$$
t_{AB}=\prod_i\left(  t_{\alpha_1}^{n_1}\circ \cdots \circ t_{\alpha_q}^{n_q}\circ   
t_{\beta_1}^{-m_1}\circ \cdots \circ t_{\beta_r}^{-m_r}   \right).
$$
\noindent Indeed, by construction we notice that $\tau_{\hat{\alpha}_i}$ and $\tau_{\hat{\beta}_j}$ are lifts of $t_{\alpha_i}$ and $t_{\beta_j}$ respectively. One obtains 
$$
\varrho\circ \tau_{\hat{\alpha}_i}=t_{\alpha_i}\circ \varrho\ \ \mbox{and}\ \ \varrho\circ \tau_{\hat{\beta}_j}=t_{\alpha_j}\circ \varrho, 
$$
\noindent which further implies that $\varrho\circ \hat{T}_{AB}=t_{AB}\circ \varrho$, as claimed.  We proved the most part of:

\noindent {\bf Lemma 6.1. } {\em Let $\Delta$ be a fixed fundamental region of $G$. Let $z\in S\backslash \{A, B\}$ and let $\hat{z}=\{\varrho^{-1}(z)\}\cap \Delta$. Then for $1\leq i\leq q$ and $1\leq j\leq r$, the geodesics $\hat{\alpha}_i$ and $\hat{\beta}_j$ in ${\Bbb{H}}$, and hence the maximal elements $D_{\hat{\alpha}_i}$ and $D_{\hat{\beta}_j}$, can be chosen to satisfy: 

\noindent $(1)$ $\varrho(\hat{\alpha}_i)=\alpha_i$, $\varrho(\hat{\beta}_j)=\beta_j$;\\
$(2)$ $\hat{\alpha}_i\cap \Delta$ and $\hat{\beta}_j\cap \Delta$ are not empty, and \\
$(3)$ the map $\hat{T}_{AB}$ defined as $(\ref{MAP})$ is a lift of $t_{AB}$ and fixes the point $\hat{z}$. }

\noindent {\em Proof. } The map $\hat{T}_{AB}$ projects to $w$ defined in the introduction, which has the property that $w(z)=z$. Thus $\hat{\alpha}_i$ and $\hat{\beta}_j$ can be chosen so that (1) and (2) are satisfied. Moreover, once $\hat{\alpha}_i$ and $\hat{\beta}_j$ are selected, $D_{\hat{\alpha}_i}$ and $D_{\hat{\beta}_j}$ can also be selected so that (\ref{ALPHA}) and (\ref{BETA}) are satisfied. 

Let $z$ lie in one component $R$ of (\ref{COL}). $R$ is either $P_i$ or $Q_j$. Since $\varrho:{\Bbb{H}}\rightarrow S$ is a local homeomorphism, there is a non-empty subset $\Sigma$ of $\Delta$ so that $\hat{z}\in \Sigma$ and $\varrho|_{\Sigma}(\Sigma)=R$. As we remarked earlier, there are perhaps more than one choices of each geodesic $\hat{\alpha}_i$ that meet $\Delta$ so that $\varrho(\hat{\alpha}_i)=\alpha_i$. In any case, there are only finitely many maximal elements of $\tau_{\hat{\alpha}_i}$ and $\tau_{\hat{\beta}_j}$ that intersect with $\Delta$. $\Sigma$ can be obtained from the fundamental region $\Delta$ with the removal of all such (finitely many) maximal elements of $\tau_{\hat{\alpha}_i}$ and $\tau_{\hat{\beta}_j}$ for $1\leq i\leq q$ and $1\leq j\leq r$.  Now from the construction we see that the restriction of $\hat{T}_{AB}$ to $\Sigma$ is the identity. In particular, we obtain $\hat{T}_{AB}(\hat{z})=\hat{z}$. This proves Lemma 6.1.  \ \ \ \ \ \ \ $\Box$

In the case of  $R=Q_j$ for some $j$, $1\leq j\leq r$, the set $\Sigma$ touches $\partial {\Bbb{H}}$ at the fixed point of a parabolic element of $G$ corresponding to the puncture $z_j$ of $Q_j$. The following lemma handles this case. 

\noindent {\bf Lemma 6.2. } {\em Under the above condition, the map $\hat{T}_{AB}$ fixes both $\hat{z}$ and the fixed point of a parabolic element of $G$. }

\noindent {\em Proof. } From Lemma 6.1, the map $\hat{T}_{AB}$ fixes $\hat{z}$, where $\hat{z}\in \Sigma$. Note that  the boundary of $\Sigma$ consists of portions of some translations of $\hat{\alpha}_i$ and $\hat{\beta}_j$. Let $z_j$ be the puncture of $Q_j$. we can draw a path $\gamma$ in $Q_j$ that connects from $z$ to $z_j$ without intersecting any boundary components of $Q_j$. In particular, $\gamma$ is disjoint from any element in $A$ or $B$.    

Now we can lift the path $\gamma$ to a path $\hat{\gamma}$ in ${\Bbb{H}}$ that connects from $\hat{z}$ to a parabolic vertex $v_j$ of $\Delta$ (corresponding the puncture $z_j$). Since $\gamma$ does not intersect with $\{A, B\}$, $\hat{\gamma}$ avoids all maximal elements  of $\tau_{\hat{\alpha}_i}$ and $\tau_{\hat{\beta}_j}$ for $1\leq i\leq q$ and $1\leq j\leq r$. But since  $\hat{T}_{AB}$ fixes $\hat{z}$ as well as any other points in $\hat{\gamma}$, by continuity, we conclude that $\hat{T}_{AB}$ fixes $v_j$, as asserted.\ \ \ \ \ $\Box$

\noindent {\em Remark $6.1$. } If $R=P_i$ for some $i$, $1\leq i\leq q$, then the region $\Sigma$ in Lemma 6.1 stays away from the boundary $\partial {\Bbb{H}}$. In this case, the argument of Lemma 6.2 fails to conclude that there is a parabolic fixed point for the map $\hat{T}_{AB}$. As a matter of fact, the Theorem of Bers \cite{Bers2} can be used to prove that $\hat{T}_{AB}$ has no any parabolic fixed point on $\partial {\Bbb{H}}$. 

%%%%%%%%%%%%%%%%%%%%
  
\section{Proof of Theorem 2}
\setcounter{equation}{0}

We assume that $S$ is non compact. recall that $Q_1,\ldots, Q_r$ obtained from (1.1) are all possible once punctured disk components of $S\backslash \{A,B\}$. 

By a theorem of Penner \cite{Pen3}, a word $w$ generated by Dehn twists along elements of $A$ and inverses of Dehn twists along elements in $B$ represents a pseudo-Anosov mapping class. So $w$ is isotopic to a pseudo-Anosov map $\omega$. Let $\chi\in \mbox{Mod}_{S}$ be induced by $\omega$, $\chi$ is hyperbolic in the sense of Bers \cite{Bers2}. It follows from Theorem 5 of Bers \cite{Bers2}, there is a Teichm\"{u}ller  geodesic $\mathcal{L}$ in $T(S)$ so that $\chi(\mathcal{L})=\mathcal{L}$. 

Assume that $S\in \mathcal{L}$. Then $\omega:S\rightarrow S$ is an absolutely extremal Teichm\"{u}ller mapping. Let $z_0\in \bar{S}$ be a zero of $\phi$ (there is nothing to prove if $\phi$ has no zeros). Since the number of zeros of $\phi$ is finite, and $\omega$ send a zero to a zero, by taking a suitable power if necessary, one may assume without loss of generality that $\omega$ fixes all zeros of $\phi$. Note that some zeros could be punctures of $S$. 

Suppose that $z_0\in Q_1$ is a non-puncture zero of $\phi$.
Let $z_1$ denote the puncture of $Q_1$. Then $z_0\neq z_1$. For any point $\hat{z}_0\in \varrho^{-1}(z_0)$, we can choose a fundamental region $\Delta$ of $G$ so that $\hat{z}_0\in \Delta$. Since $Q_1$ contains a puncture, there is a parabolic vertex $v_1$ of $\Delta$ in $\partial {\Bbb{H}}$ that corresponds to the puncture $z_1$. 

From Lemma 6.1 and Lemma 6.2, $w$ can be lifted to $\hat{T}_{AB}$ that fixes $\hat{z}_0$ and $v_1$. Since $\omega$ is isotopic to $w$, $\omega$ can also be lifted to $\hat{\omega}$ so that $\hat{\omega}|_{\partial {\Bbb{H}}}=\hat{T}_{AB}|_{\partial {\Bbb{H}}}$. It follows that $[\hat{\omega}]=[\hat{T}_{AB}]$ and $\hat{\omega}$ fixes both $\hat{z}_0$ and $v_1$ as well. 

By Lemmas 5.2, $[\hat{\omega}]^*\in \mbox{Mod}_S^a$ is not pseudo-Anosov. More precisely, let $[\hat{\omega}]^*$ be induced by $W:\dot{S}\rightarrow \dot{S}$. By Lemma 5.1, $W$ fixes the boundary of a twice punctured disk $\Omega$ enclosing $a$ so that $W|_{\Omega}=\mbox{id}$ and $W|_{\dot{S}\backslash D}$ is pseudo-Anosov. In particular, $W$ is reducible. Set $\theta=[\hat{\omega}]=[\hat{T}_{AB}]$. We claim that $\theta$ satisfies all conditions of Theorem 2. 

Indeed, (1) is clear. By definition, $\theta^*$ is not  pseudo-Anosov, which satisfies (3). Finally, since $\hat{T}_{AB}$ fixes $\hat{z}_0$, by Lemma 3.2, $\theta$ keeps the real one dimensional submanifold 
$\hat{\mathcal{L}}\subset F(S)$ (defined in (3.1)) invariant. So (2) holds. This completes the proof of Theorem 2, and hence of Theorem 1.  \ \ \ $\Box$

\end{document}